\documentclass[12pt]{article}
\usepackage{amsmath,amssymb,amsthm}
\usepackage{hyperref}
\usepackage{datetime}
\usepackage{units}
\usepackage{color}
\usepackage[T1]{fontenc}
\usepackage[utf8]{inputenc}
\usepackage{authblk}
\usepackage{bm,latexsym,mathrsfs,enumerate}
\usepackage{color}

\title{ A novel approach to the discovery of ternary BBP-type formulas for polylogarithm constants\thanks{%
MSC 2010: 11Y60, 30B99}}
\author[]{Kunle Adegoke\thanks{adegoke00@gmail.com }}
\affil{Department of Physics and Engineering Physics, \mbox{Obafemi Awolowo University, Ile-Ife, 220005 Nigeria}}

\theoremstyle{plain}
\numberwithin{equation}{section}
\begin{document}
\date{}
\maketitle
%\today, \currenttime
\begin{abstract}
\noindent Using a clear and straightforward approach, we prove new ternary (base 3) digit extraction BBP-type formulas for polylogarithm constants. Some known results are also rediscovered in a more direct and elegant manner. A previously unproved degree~4 ternary formula is also proved. Finally, a couple of ternary zero relations are established, which prove two known but hitherto unproved formulas.
\end{abstract}
\tableofcontents

\section{Introduction}	

BBP-type formulas are formulas of the form
\[
\alpha = \sum_{k=0}^\infty  1/b^k \sum_{j=1}^n a_j / (k n + j)^s
\]
where $s$, $b$, $n$ and $a_j$ are integers, and $\alpha$ is some constant. Formulas of this type were first introduced in a 1996 paper~\cite{BBP97}, where a formula of this type for $\pi$ was given.  Such formulas have the remarkable property that they permit one to calculate \mbox{base-$b$} digits of the constant $\alpha$ beginning at an arbitrary starting position, by means of a simple algorithm that requires almost no memory and (depending on how many digits are required) without the need for multiple-precision arithmetic software~\cite{bailey09}.  Such formulas also have intriguing connections to the \mbox{age-old} problem of understanding why the digits of various transcendental constants appear ``normal'' -- each string of \mbox{$m$-long} digits appears, in the limit, with frequency $1/b^m$~\cite{bailey09,bailey01,borwein02,chamberland03}.

\bigskip

While many binary BBP-type formulas are now known, only relatively few ternary (base-3) BBP-type formulas have been discovered. This present paper is concerned with the symbolic (that is, non-computer-search-based) discovery of ternary (base-3) BBP-type formulas for polylogarithm constants. The methods used here aim to complement the experimental approaches that have dominated the area. Through fundamental methods, a wide range of interesting formulas will be obtained. In most cases, the procedure for obtaining the ternary formulas shall consist mainly of evaluating a polylogarithm functional equation at indicated coordinates and noting the following identities for the real and imaginary parts of the polylogarithm function:
\begin{equation}\label{equ.urj5bd8}
\begin{split}
{\mathop{\rm Re\;}\nolimits} {\rm Li}_s \left[ {pe^{ix} } \right] &= \sum\limits_{k = 1}^\infty  {\frac{{p^k \cos kx}}{{k^s }}}\,,\\
{\mathop{\rm Im\;}\nolimits} {\rm Li}_s \left[ {pe^{ix} } \right] &= \sum\limits_{k = 1}^\infty  {\frac{{p^k \sin kx}}{{k^s }}}\,,
\end{split}
\end{equation}
for $p\in [0,1]$, $x\in\mathbb{R}$ and $s\in\mathbb{Z^+}$. In the above equations,  {\rm Li} is the notation for the polylogarithm function defined by
\[
{\rm Li}_s [z] = \sum\limits_{k = 1}^\infty  {\frac{{z^k }}{{k^s }}},\quad |z|\le 1\,. 
\]

When $p=1$ we have
\[
\begin{split}
{\rm Li}_{2n} [e^{ix} ]& = {\rm Gl}_{2n} (x) + i{\rm Cl}_{2n} (x)\\
{\rm Li}_{2n + 1} [e^{ix} ]& = {\rm Cl}_{2n + 1} (x) + i{\rm Gl}_{2n + 1} (x)\,,
\end{split}
\]
where ${\rm Gl}$ and ${\rm Cl}$ are Clausen sums~\cite{lewin81} defined, for $n\in\mathbb{Z^+}$ by

\[
\begin{split}
{\rm Cl}_{2n} (x)& = \sum\limits_{k = 1}^\infty  {\frac{{\sin kx}}{{k^{2n} }}},\quad {\rm Cl}_{2n + 1} (x) = \sum\limits_{k = 1}^\infty  {\frac{{\cos kx}}{{k^{2n + 1} }}}\\
{\rm Gl}_{2n} (x) &= \sum\limits_{k = 1}^\infty  {\frac{{\cos kx}}{{k^{2n} }}},\quad {\rm Gl}_{2n + 1} (x) = \sum\limits_{k = 1}^\infty  {\frac{{\sin kx}}{{k^{2n + 1} }}}\,.
\end{split}
\]

We shall find the following formulas useful:

\[
\begin{split}
{\rm Gl}_{2n} (x) &= ( - 1)^{1 + [n/2]} 2^{n - 1} \pi ^n {\rm B}_n (x/2\pi )/n!\\
\frac{1}{{m^{n - 1} }}{\rm Cl}_n (mx) &= \sum\limits_{r = 0}^{m - 1} {{\rm Cl}_{n} (x + 2\pi r/m)}\,.
\end{split}
\]
Here $[n/2]$ denotes the integer part of $n/2$ and $ {\rm B}_n$ are the Bernoulli polynomials defined by

\[
{\frac {t{{\rm e}^{xt}}}{{{\rm e}^{t}}-1}}=\sum _{n=0}^{\infty }{
\frac {{\rm B}_n (x) {t}^{n}}{n!}}\,.
\]

\section{Degree~$1$ Ternary BBP-type Formulas}
In reference~\cite{adegokeBBP2}, several degree~$1$ BBP-type formulas in general bases are proven. In many of the formulas, ternary formulas may be readily obtained by writing the base in each case as a power of $3$. 

\bigskip

Here we now present a couple of interesting degree~$1$ ternary formulas.

\bigskip

The following identities are easily verified:
\[
\begin{split}
{\rm Li}_1 \left[ {\frac{1}{{\sqrt 3 }}\exp \left( {\frac{{i\pi }}{6}} \right)} \right] &= \frac{1}{2}\ln 3 + \frac{{i\pi }}{6}\\
\text{and}\\
{\rm Li}_1 \left[ {\frac{1}{{\sqrt 3 }}\exp \left( {\frac{{i\pi }}{2}} \right)} \right] &= \frac{1}{2}\ln 3 - \ln 2 + \frac{{i\pi }}{6}\,.
\end{split}
\]

We therefore have the formulas:

\begin{equation}\label{equ.gtkbg3t}
\ln 2 = {\mathop{\rm Re\,}\nolimits} {\rm Li}_1 \left[ {\frac{1}{{\sqrt 3 }}\exp \left( {\frac{{i\pi }}{6}} \right)} \right] - {\mathop{\rm Re\,}\nolimits} {\rm Li}_1 \left[ {\frac{1}{{\sqrt 3 }}\exp \left( {\frac{{i\pi }}{2}} \right)} \right]\,,
\end{equation}

\begin{equation}\label{equ.zh8e577}
\ln 3=2 {\mathop{\rm Re\,}\nolimits} {\rm Li}_1 \left[ {\frac{1}{{\sqrt 3 }}\exp \left( {\frac{{i\pi }}{6}} \right)} \right]
\end{equation}
and

\begin{equation}\label{equ.xwkmwte}
\pi = 6\,{\mathop{\rm Im\,}\nolimits} {\rm Li}_1 \left[ {\frac{1}{{\sqrt 3 }}\exp \left( {\frac{{i\pi }}{6}} \right)} \right] = 6\,{\mathop{\rm Im\,}\nolimits} {\rm Li}_1 \left[ {\frac{1}{{\sqrt 3 }}\exp \left( {\frac{{i\pi }}{2}} \right)} \right]\,.
\end{equation}

It is also straightforward to verify that:
\begin{equation}\label{equ.c8cayj2}
\begin{split}
\ln 2 &= {\rm Li}_1 \left[ {\frac{1}{3}} \right] - {\rm Li}_1 \left[ { - \frac{1}{3}} \right]\\
\text{and}\\ 
\ln 3 &= 2{\rm Li}_1 \left[ {\frac{1}{3}} \right] - {\rm Li}_1 \left[ { - \frac{1}{3}} \right]\,. 
\end{split}
\end{equation}
Based on the above identities, we are now ready to derive explicit BBP-type formulas for $\ln 2$, $\ln 3$ and $\pi$.
\subsection{Ternary formulas for $\ln 2$}

Using the first equality of~\eqref{equ.urj5bd8}, we note that

\begin{equation}\label{equ.tmvvkod}
\begin{split}
{\mathop{\rm Re\,}\nolimits} {\rm Li}_{\rm 1} \left[ {\frac{1}{{\sqrt 3 }}\exp \left( {\frac{{i\pi }}{6}} \right)} \right] &= \sum\limits_{k = 1}^\infty  {\left( {\frac{1}{{\sqrt 3 }}} \right)^k \frac{{\cos (k\pi /6)}}{k}}\\
&= \frac{1}{{2 \cdot 3^6 }}\sum\limits_{k = 0}^\infty  {\frac{1}{{3^{6k} }}\left[ {\frac{{3^6 }}{{12k + 1}} + \frac{{3^5 }}{{12k + 2}}} \right.}\\
&\qquad- \frac{{3^4 }}{{12k + 4}} - \frac{{3^4 }}{{12k + 5}} - \frac{{2 \cdot 3^3 }}{{12k + 6}} - \frac{{3^3 }}{{12k + 7}}\\
&\qquad\qquad\left. { - \frac{{3^2 }}{{12k + 8}} + \frac{3}{{12k + 10}} + \frac{3}{{12k + 11}} + \frac{2}{{12k + 12}}} \right]
\end{split}
\end{equation}
and

\begin{equation}\label{equ.hlnz76f}
\begin{split}
{\mathop{\rm Re\;}\nolimits} {\rm Li}_{\rm 1} \left[ {\frac{1}{{\sqrt 3 }}\exp \left( {\frac{{i\pi }}{2}} \right)} \right] &= \frac{1}{{3^6 }}\sum\limits_{k = 0}^\infty  {\frac{1}{{3^{6k} }}\left[ { - \frac{{3^5 }}{{12k + 2}} + \frac{{3^4 }}{{12k + 4}}} \right.}\\
&\qquad- \frac{{3^3 }}{{12k + 6}}\left. { + \frac{{3^2 }}{{12k + 8}} - \frac{3}{{12k + 10}} + \frac{2}{{12k + 12}}} \right]\,.
\end{split}
\end{equation}
Subtracting~\eqref{equ.hlnz76f} from~\eqref{equ.tmvvkod} in accordance with~\eqref{equ.gtkbg3t}, we obtain the following ternary BBP-type formula for $\ln 2$:
\begin{equation}\label{equ.m5b8en9}
\begin{split}
\ln 2 &= \frac{1}{{2 \cdot 3^5 }}\sum\limits_{k = 0}^\infty  {\frac{1}{{3^{6k} }}\left[ {\frac{{3^5 }}{{12k + 1}} + \frac{{3^5 }}{{12k + 2}}} \right.}\\
&\quad- \frac{{3^4 }}{{12k + 4}} - \frac{{3^3 }}{{12k + 5}} - \frac{{3^2 }}{{12k + 7}} - \frac{{3^2 }}{{12k + 8}}\\
&\qquad\left. { + \frac{3}{{12k + 10}} + \frac{1}{{12k + 11}}} \right]\,.
\end{split}
\end{equation}

Note that an alternating version of~\eqref{equ.m5b8en9}, using the same scheme, is

\[
\ln 2 = \frac{1}{{18}}\sum\limits_{k = 0}^\infty  {\left( { - \frac{1}{{27}}} \right)^k \left[ {\frac{9}{{6k + 1}} + \frac{9}{{6k + 2}} - \frac{3}{{6k + 4}} - \frac{1}{{6k + 5}}} \right]}\,. 
\]
From the first equality of~\eqref{equ.c8cayj2}, we can obtain yet another ternary formula for $\ln 2$. We first note that
\begin{equation}\label{equ.fv07456}
{\rm Li}_1 \left[ {\frac{1}{3}} \right] = \frac{1}{9}\sum\limits_{k = 0}^\infty  {\frac{1}{{9^k }}\left[ {\frac{3}{{2k + 1}} + \frac{1}{{2k + 2}}} \right]} 
\end{equation}
and

\begin{equation}\label{equ.rm10gx0}
{\rm Li}_1 \left[ { - \frac{1}{3}} \right] = \frac{1}{9}\sum\limits_{k = 0}^\infty  {\frac{1}{{9^k }}\left[ { - \frac{3}{{2k + 1}} + \frac{1}{{2k + 2}}} \right]}\,. 
\end{equation}
Subtracting~\eqref{equ.rm10gx0} from~\eqref{equ.fv07456} in accordance with the first equality of~\eqref{equ.c8cayj2}, we obtain
\begin{equation}\label{equ.vnbfew8}
\ln 2 = \frac{2}{3}\sum\limits_{k = 0}^\infty  {\frac{1}{{9^k }}\left[ {\frac{1}{{2k + 1}}} \right]}\,, 
\end{equation}
which is listed in the BBP Compendium as formula~(64).

\subsection{Ternary formulas for $\ln 3$}

From~\eqref{equ.tmvvkod} and~\eqref{equ.zh8e577}, we obtain the following ternary BBP-type formula for $\ln 3$:

\begin{equation}\label{equ.awtepfw}
\begin{split}
\ln 3 &= \frac{1}{{3^6 }}\sum\limits_{k = 0}^\infty  {\frac{1}{{3^{6k} }}\left[ {\frac{{3^6 }}{{12k + 1}} + \frac{{3^5 }}{{12k + 2}}} \right.}\\
&\quad- \frac{{3^4 }}{{12k + 4}} - \frac{{3^4 }}{{12k + 5}} - \frac{{2 \cdot 3^3 }}{{12k + 6}} - \frac{{3^3 }}{{12k + 7}}\\
&\qquad\left. { - \frac{{3^2 }}{{12k + 8}} + \frac{3}{{12k + 10}} + \frac{3}{{12k + 11}} + \frac{2}{{12k + 12}}} \right]\,.
\end{split}
\end{equation}

An alternating version of the above formula is

\[
\ln 3 = \frac{1}{{27}}\sum\limits_{k = 0}^\infty  {\left( { - \frac{1}{{27}}} \right)^k \left[ {\frac{{27}}{{6k + 1}} + \frac{9}{{6k + 2}} - \frac{3}{{6k + 4}} - \frac{3}{{6k + 5}} - \frac{2}{{6k + 6}}} \right]}\,. 
\]

Combining~\eqref{equ.fv07456} and~\eqref{equ.rm10gx0} according to the second equality of~\eqref{equ.c8cayj2} we obtain another ternary BBP-type formula for $\ln 3$ as :
\begin{equation}\label{equ.rbqo59l}
\ln 3 = \frac{1}{9}\sum\limits_{k = 0}^\infty  {\frac{1}{{9^k }}\left[ {\frac{9}{{2k + 1}} + \frac{1}{{2k + 2}}} \right]}\,,
\end{equation}
which is Formula~(67) of the Compendium.

\subsection{Ternary formulas for $\pi\sqrt 3$}

From~\eqref{equ.xwkmwte}, we immediately obtain the ternary BBP-type formulas
\begin{equation}\label{equ.vd5d381}
\begin{split}
\pi\sqrt 3  &= \frac{{1 }}{{3^4 }}\sum\limits_{k = 0}^\infty  {\frac{1}{{3^{6k} }}\left[ {\frac{{3^5 }}{{12k + 1}} + \frac{{3^5 }}{{12k + 2}} + \frac{{2 \cdot 3^4 }}{{12k + 3}}} \right.}\\
&\quad+ \frac{{3^4 }}{{12k + 4}} + \frac{{3^3 }}{{12k + 5}} - \frac{{3^2 }}{{12k + 7}} - \frac{{3^2 }}{{12k + 8}}\\
&\qquad\left. { - \frac{{2 \cdot 3}}{{12k + 9}} - \frac{3}{{12k + 10}} - \frac{1}{{12k + 11}}} \right]
\end{split}
\end{equation}

and

\begin{equation}\label{equ.whhwzp3}
\begin{split}
\pi\sqrt 3  &= \frac{{2 }}{{3^4 }}\sum\limits_{k = 0}^\infty  {\frac{1}{{3^{6k} }}\left[ {\frac{{3^5 }}{{12k + 1}} - \frac{{3^4 }}{{12k + 3}} + \frac{{3^3 }}{{12k + 5}}} \right.}\\
&\qquad- \frac{{3^2 }}{{12k + 7}}\left. { + \frac{3}{{12k + 9}} - \frac{1}{{12k + 11}}} \right]\,.
\end{split}
\end{equation}

An alternating version of~\eqref{equ.whhwzp3} is
\begin{equation}
\pi\sqrt 3  = 6 \sum\limits_{k = 0}^\infty  {\left( { - \frac{1}{3}} \right)^k \left[ {\frac{1}{{2k + 1}}} \right]}\,. 
\end{equation}

\subsection{Ternary Zero Relations}
Identity~\eqref{equ.vnbfew8} may be rewritten in base~$3^6$, length~$12$ as
\begin{equation}\label{equ.wrwrv1q}
\ln 2 = \frac{4}{{3^5 }}\sum\limits_{k = 0}^\infty  {\frac{1}{{3^{6k} }}\left[ {\frac{{3^4 }}{{12k + 2}} + \frac{{3^2 }}{{12k + 6}} + \frac{{1 }}{{12k + 10}}} \right]}\,.
\end{equation}

Subtracting~\eqref{equ.wrwrv1q} from~\eqref{equ.m5b8en9}, we obtain the following ternary zero relation:

\begin{equation}\label{equ.f51s4o9}
\begin{split}
&0 = \sum\limits_{k = 0}^\infty  {\frac{1}{{3^{6k} }}\left[ {\frac{{3^5}}{{12k + 1}} - \frac{{5\cdot 3^4}}{{12k + 2}} - \frac{{3^4}}{{12k + 4}}} \right.}\\
&\qquad\qquad-\frac{{3^3}}{{12k + 5}} - \frac{{2^3\cdot 3^2}}{{12k + 6}} - \frac{3^2}{{12k + 7}} - \frac{3^2}{{12k + 8}}\\
&\qquad\qquad\qquad\left. { - \frac{5}{{12k + 10}} + \frac{1}{{12k + 11}}} \right]\,.
\end{split}
\end{equation}
Let $R127$ denote the right hand side of Compendium formula~127 and $R128$ the right hand side of Compendium formula~128. We note that 

\begin{equation}\label{equ.k2e5m43}
{\rm equation}~\eqref{equ.f51s4o9}\equiv R128-R127=0\,. 
\end{equation}

Subtracting~\eqref{equ.vd5d381} from~\eqref{equ.whhwzp3}, we obtain the zero relation:

\begin{equation}\label{equ.bi0kz7u}
\begin{split}
0 &= \sum\limits_{k = 0}^\infty  {\frac{1}{{3^{6k} }}\left[ {\frac{{243}}{{12k + 1}} - \frac{{243}}{{12k + 2}} - \frac{{324}}{{12k + 3}}} \right.}\\
&\qquad- \frac{{81}}{{12k + 4}} + \frac{{27}}{{12k + 5}} - \frac{9}{{12k + 7}} + \frac{9}{{12k + 8}}\\
&\qquad\qquad+ \frac{{12}}{{12k + 9}}\left. { + \frac{3}{{12k + 10}} - \frac{1}{{12k + 11}}} \right]\,.
\end{split}
\end{equation}

Note also that
\begin{equation}\label{equ.vakij8r}
{\rm equation}~\eqref{equ.bi0kz7u}\equiv R128 + R127=0 
\end{equation}

Equations~\eqref{equ.k2e5m43} and \eqref{equ.vakij8r} therefore establish that

\[
\begin{split}
R127&=0\\ 
\mbox{and}\\
R128 &=0\,. 
\end{split}
\]
Thus the hitherto unproved formulas (127) and (128) in the BBP Compendium are now proved.

\bigskip

\section{Degree~$2$ Ternary BBP-type Formulas}

The dilogarithm reflection formula (\mbox{equation~A.2.1.7} of~\cite{lewin81}) is 

\[
\frac{{\pi ^2 }}{6} - \ln x\ln (1 - x)={\rm Li}_2 [x] + {\rm Li}_2 [1 - x]\,.
\]
Putting $x=-\exp(i\pi/3)$ in the above identity and taking real and imaginary parts we find
\begin{equation}\label{equ.p12ua6j}
\frac{{5\pi ^2 }}{{72}} - \frac{1}{8}\ln ^2 3={\mathop{\rm Re\,}\nolimits} {\rm Li}_2 \left[ {\frac{1}{{\sqrt 3 }}\exp \left( {\frac{{i\pi }}{6}} \right)} \right]
\end{equation}

and

\begin{equation}\label{equ.og6rqzy}
\frac{2}{3}{\rm Cl}_2 \left( {\frac{\pi }{3}} \right) - \frac{{\pi \ln 3}}{{12}}={\mathop{\rm Im\,}\nolimits} {\rm Li}_2 \left[ {\frac{1}{{\sqrt 3 }}\exp \left( {\frac{{i\pi }}{6}} \right)} \right]\,.
\end{equation}

A two-variable functional equation for dilogarithms, due to Kummer (\mbox{equation~A.2.1.19} of~\cite{lewin81}) is

\begin{equation}\label{equ.zbzsvtl}
\begin{split}
{\rm Li}_2 \left[ {\frac{{x(1 - y)^2 }}{{y(1 - x)^2 }}} \right] &= {\rm Li}_2 \left[ { - \frac{{x(1 - y)}}{{(1 - x)}}} \right] + {\rm Li}_2 \left[ { - \frac{{(1 - y)}}{{y(1 - x)}}} \right]\\
&\qquad+ {\rm Li}_2 \left[ {\frac{x}{y}\frac{{(1 - y)}}{{(1 - x)}}} \right] + {\rm Li}_2 \left[ {\frac{{1 - y}}{{1 - x}}} \right] + \frac{1}{2}\ln ^2 y\,.
\end{split}
\end{equation}

Choosing $x=-\exp(i\pi/3)$ and $y=\exp(i\pi/3)$ in~\eqref{equ.zbzsvtl} gives

\begin{equation}\label{equ.q21efos}
\frac{{\pi ^2 }}{{12}} - \frac{{\ln ^2 3}}{4}={\rm Li}_2 \left[ {\frac{1}{3}} \right] - \frac{1}{2}{\rm Li}_2 \left[ { - \frac{1}{3}} \right]\,.
\end{equation}

Note that the choice of $x=-1$ and $y=1/3$ gives the same result.

\bigskip

Another two-variable functional equation for dilogarithms, due to Abel (\mbox{equation~A.2.1.16} of~\cite{lewin81}) is

\begin{equation}\label{equ.u7rchcl}
\begin{split}
{\rm Li}_2 \left[ {\frac{x}{{1 - x}} \cdot \frac{y}{{1 - y}}} \right] &= {\rm Li}_2 \left[ {\frac{x}{{(1 - y)}}} \right] + {\rm Li}_2 \left[ {\frac{y}{{(1 - x)}}} \right]\\
&\qquad- {\rm Li}_2 \left[ x \right] - {\rm Li}_2 \left[ y \right] - \ln (1 - x)\ln (1 - y)\,.
\end{split}
\end{equation}

Choosing $x=-\exp(i\pi/3)$ and $y=\exp(-i\pi/3)$ in~\eqref{equ.u7rchcl} and taking imaginary parts, we obtain

\begin{equation}\label{equ.qiq0qug}
\frac{5}{2}{\rm Cl}_2 \left( {\frac{\pi }{3}} \right) - \frac{{\pi \ln 3}}{4} = 3{\mathop{\rm Im\,}\nolimits} {\rm Li}_2 \left[ {\frac{1}{{\sqrt 3 }}\exp \left( {\frac{{i\pi }}{2}} \right)} \right]\,.
\end{equation}

\subsection{Ternary Formula for $\pi^2$}

Solving~\eqref{equ.p12ua6j} and~\eqref{equ.q21efos} for $\pi^2$, we obtain
\begin{equation}\label{equ.kz652nt}
\pi ^2  = 36{\mathop{\rm Re\,}\nolimits} {\rm Li}_2 \left[ {\frac{1}{{\sqrt 3 }}\exp \left( {\frac{{i\pi }}{6}} \right)} \right] - 18{\rm Li}_2 \left[ {\frac{1}{3}} \right] + 9{\rm Li}_2 \left[ { - \frac{1}{3}} \right]\,.
\end{equation}

Writing

\begin{equation}\label{equ.rrx9k6k}
\begin{split}
{\mathop{\rm Re\,}\nolimits} {\rm Li}_2 \left[ {\frac{1}{{\sqrt 3 }}\exp \left( {\frac{{i\pi }}{6}} \right)} \right] &= \frac{1}{{2 \cdot 3^6 }}\sum\limits_{k = 0}^\infty  {\frac{1}{{3^{6k} }}\left[ {\frac{{3^6 }}{{(12k + 1)^2 }} + \frac{{3^5 }}{{(12k + 2)^2 }}} \right.}\\
&\qquad- \frac{{3^4 }}{{(12k + 4)^2 }} - \frac{{3^4 }}{{(12k + 5)^2 }} - \frac{{3^3 }}{{(12k + 6)^2 }}\\
&\qquad\qquad- \frac{{3^3 }}{{(12k + 7)^2 }} - \frac{{3^2 }}{{(12k + 8)^2 }}\\
&\qquad\qquad\qquad\left. { + \frac{3}{{(12k + 10)^2 }} + \frac{2}{{(12k + 12)^2 }}} \right]
\end{split}
\end{equation}

and

\begin{equation}\label{equ.ksebql9}
\begin{split}
{\rm Li}_2 \left[ { \pm \frac{1}{3}} \right]& = \frac{{2^2 }}{{3^6 }}\sum\limits_{k = 0}^\infty  {\frac{1}{{3^{6k} }}\left[ { \pm \frac{{3^5 }}{{(12k + 2)^2 }} + \frac{{3^4 }}{{(12k + 4)^2 }}} \right.} \\
&\quad\pm \frac{{3^3 }}{{(12k + 6)^2 }} + \frac{{3^2 }}{{(12k + 8)^2 }}\\
&\qquad\left. { \pm \frac{3}{{(12k + 10)^2 }} + \frac{1}{{(12k + 12)^2 }}} \right]\,,
\end{split}
\end{equation}

and combining them according to~\eqref{equ.kz652nt} we establish a ternary BBP-type formula for $\pi^2$:
\begin{equation}\label{equ.rlkmbcg}
\begin{split}
\pi ^2 & = \frac{2}{{27}}\sum\limits_{k = 0}^\infty  {\frac{1}{{3^{6k} }}\left[ {\frac{{3^5 }}{{(12k + 1)^2 }} - \frac{{5 \cdot 3^4 }}{{(12k + 2)^2 }}} \right.} \\
&\quad- \frac{{3^4 }}{{(12k + 4)^2 }} - \frac{{3^3 }}{{(12k + 5)^2 }} - \frac{{2^3  \cdot 3^2 }}{{(12k + 6)^2 }}\\
&\qquad\qquad- \frac{{3^2 }}{{(12k + 7)^2 }} - \frac{{3^2 }}{{(12k + 8)^2 }}\\
&\qquad\qquad\qquad\left. { - \frac{5}{{(12k + 10)^2 }} + \frac{1}{{(12k + 11)^2 }}} \right]\,.
\end{split}
\end{equation}

Incidentally,~\eqref{equ.rlkmbcg} is Formula~(73) of the Compendium.

\subsection{Ternary Formula for $\ln^2 3$}

Solving~\eqref{equ.p12ua6j} and~\eqref{equ.q21efos} for $\ln^2 3$, we obtain
\begin{equation}\label{equ.k0ed0v7}
 \ln^2 3  = 12\,{\mathop{\rm Re\,}\nolimits} {\rm Li}_2 \left[ {\frac{1}{{\sqrt 3 }}\exp \left( {\frac{{i\pi }}{6}} \right)} \right] - 10\,{\rm Li}_2 \left[ {\frac{1}{3}} \right] + 5\,{\rm Li}_2 \left[ { - \frac{1}{3}} \right]\,.
\end{equation}

Using~\eqref{equ.rrx9k6k} and \eqref{equ.ksebql9} above in~\eqref{equ.k0ed0v7}, we obtain the ternary BBP-type formula for $\ln^2 3$ as

\begin{equation}\label{equ.cwx7fio}
\begin{split}
\ln^2 3 & = \frac{1}{{3^6}}\sum\limits_{k = 0}^\infty  {\frac{1}{{3^{6k} }}\left[ {\frac{{2\cdot 3^7 }}{{(12k + 1)^2 }} - \frac{{2\cdot 3^8 }}{{(12k + 2)^2 }}} \right.} \\
&\quad- \frac{{2\cdot 13\cdot 3^4 }}{{(12k + 4)^2 }} - \frac{{2\cdot 3^5 }}{{(12k + 5)^2 }} - \frac{{2^3  \cdot 3^5 }}{{(12k + 6)^2 }}\\
&\qquad\qquad- \frac{{2\cdot 3^4 }}{{(12k + 7)^2 }} - \frac{{2\cdot 13\cdot 3^2 }}{{(12k + 8)^2 }}\\
&\qquad\qquad\qquad\left. { - \frac{2\cdot 3^4}{{(12k + 10)^2 }} + \frac{2\cdot 3^2}{{(12k + 11)^2 }}- \frac{8}{{(12k + 12)^2 }}} \right]\,,
\end{split}
\end{equation}

which is Formula~(74) of the Compendium.

\subsection{Ternary Formula for $\pi\sqrt 3\ln 3$}

Solving~\eqref{equ.og6rqzy} and \eqref{equ.qiq0qug} for $\pi\ln 3$ gives 

\begin{equation}\label{equ.lswk9s4}
\pi \ln 3 = 48\,{\mathop{\rm Im\,}\nolimits} {\rm Li}_2 \left[ {\frac{1}{{\sqrt 3 }}\exp \left( {\frac{{i\pi }}{2}} \right)} \right] - 60\,{\mathop{\rm Im\,}\nolimits} {\rm Li}_2 \left[ {\frac{1}{{\sqrt 3 }}\exp \left( {\frac{{i\pi }}{6}} \right)} \right]\,.
\end{equation}

Now

\begin{equation}\label{equ.x7uw9cw}
\begin{split}
{\mathop{\rm Im\,}\nolimits} {\rm Li}_2 \left[ {\frac{1}{{\sqrt 3 }}\exp \left( {\frac{{i\pi }}{6}} \right)} \right] &= \frac{{\sqrt 3 }}{{54}}\sum\limits_{k = 0}^\infty  {\left( { - \frac{1}{{27}}} \right)^k \left[ {\frac{9}{{(6k + 1)^2 }} + \frac{9}{{(6k + 2)^2 }}} \right.}\\
&\qquad\left. { + \frac{6}{{(12k + 3)^2 }} + \frac{3}{{(12k + 4)^2 }} + \frac{1}{{(6k + 5)^2 }}} \right]
\end{split}
\end{equation}

and

\begin{equation}\label{equ.mjrwr77}
{\mathop{\rm Im\,}\nolimits} {\rm Li}_2 \left[ {\frac{1}{{\sqrt 3 }}\exp \left( {\frac{{i\pi }}{2}} \right)} \right] = \frac{{\sqrt 3 }}{{27}}\sum\limits_{k = 0}^\infty  {\left( { - \frac{1}{{27}}} \right)^k \left[ {\frac{9}{{(6k + 1)^2 }} - \frac{3}{{(6k + 3)^2 }} + \frac{1}{{(6k + 5)^2 }}} \right]}
\end{equation}

Using~\eqref{equ.x7uw9cw} and \eqref{equ.mjrwr77} in~\eqref{equ.lswk9s4} leads to the ternary BBP-type formula

\[
\begin{split}
\pi \sqrt 3\ln 3 &= 2\sum\limits_{k = 0}^\infty  {\left( { - \frac{1}{{27}}} \right)^k \left[ {\frac{9}{{(6k + 1)^2 }} - \frac{{15}}{{(6k + 2)^2 }} - \frac{{18}}{{(6k + 3)^2 }}} \right.}\\
&\qquad\left. { - \frac{5}{{(6k + 4)^2 }} + \frac{1}{{(6k + 5)^2 }}} \right]\,.
\end{split}
\]

\subsection{Ternary Formula for $\sqrt 3 {\rm Cl}_2 (\pi /3)$}

Solving~\eqref{equ.og6rqzy} and \eqref{equ.qiq0qug} for ${\rm Cl}_2 (\pi /3)$ gives 

\begin{equation}\label{equ.z9u5qsv}
{\rm Cl}_2 \left( {\frac{\pi }{3}} \right) = 6\,{\mathop{\rm Im\,}\nolimits} {\rm Li}_2 \left[ {\frac{1}{{\sqrt 3 }}\exp \left( {\frac{{i\pi }}{2}} \right)} \right] - 6\,{\mathop{\rm Im\,}\nolimits} {\rm Li}_2 \left[ {\frac{1}{{\sqrt 3 }}\exp \left( {\frac{{i\pi }}{6}} \right)} \right]\,.
\end{equation}

Using~\eqref{equ.x7uw9cw} and \eqref{equ.mjrwr77} in~\eqref{equ.z9u5qsv} leads to the ternary BBP-type formula

\[
\begin{split}
\sqrt 3{\rm Cl}_2 \left( {\frac{\pi }{3}} \right) &= \frac{{1 }}{3}\sum\limits_{k = 0}^\infty  {\left( { - \frac{1}{{27}}} \right)^k \left[ {\frac{9}{{(6k + 1)^2 }} - \frac{9}{{(6k + 2)^2 }} - \frac{{12}}{{(6k + 3)^2 }}} \right.}\\
&\qquad\left. { - \frac{3}{{(6k + 4)^2 }} + \frac{1}{{(6k + 5)^2 }}} \right]\,.
\end{split}
\]

\section{Degree~$3$ Ternary BBP-type Formulas}

A functional identity for trilogarithms (\mbox{equation A.2.6.10} of~\cite{lewin81}) is
\[
\begin{split}
{\rm Li}_3 \left[ {\frac{{1 - x}}{{1 + x}}} \right]& - {\rm Li}_3 \left[ {\frac{{x - 1}}{{x + 1}}} \right] = 2\,{\rm Li}_3 \left[ {1 - x} \right] + 2\,{\rm Li}_3 \left[ {\frac{1}{{1 + x}}} \right]\\
&\qquad\qquad\qquad\qquad- \frac{1}{2}\,{\rm Li}_3 \left[ {1 - x^2 } \right] - \frac{7}{4}\,\zeta (3) \\
&\qquad\qquad\qquad\qquad\qquad+\frac{{\pi ^2 }}{6}\ln (1 + x) - \frac{1}{3}\ln ^3 (1 + x)\,.
\end{split}
\]

The use of $x=2$ in the above equation gives

\[
\frac{{13}}{6}\zeta (3) - \frac{1}{6}\pi ^2 \ln 3 + \frac{1}{6}\ln ^3 3 = 2\,{\rm Li}_3 \left[ {\frac{1}{3}} \right] - {\rm Li}_3 \left[ { - \frac{1}{3}} \right]\,.
\]

Putting $x=\exp{i\pi/3}$ in the functional equation and taking real and imaginary parts gives

\begin{equation}\label{equ.thgd6ju}
\frac{{13}}{{18}}\zeta (3) - \frac{5}{{144}}\pi ^2 \ln 3 + \frac{1}{{48}}\ln ^3 3 = {\mathop{\rm Re\,}\nolimits} {\rm Li}_3 \left[ {\frac{1}{{\sqrt 3 }}\exp \left( {\frac{{i\pi }}{6}} \right)} \right]
\end{equation}

and

\begin{equation}\label{equ.hv0t3gq}
\begin{split}
&\frac{{29}}{{1296}}\pi ^3 - \frac{1}{{48}}\pi \ln ^2 3\\
&\qquad= 4\,{\mathop{\rm Im\,}\nolimits} {\rm Li}_3 \left[ {\frac{1}{{\sqrt 3 }}\exp \left( {\frac{{i\pi }}{2}} \right)} \right] - 5\,{\mathop{\rm Im\,}\nolimits} {\rm Li}_3 \left[ {\frac{1}{{\sqrt 3 }}\exp \left( {\frac{{i\pi }}{6}} \right)} \right]\,.
\end{split}
\end{equation}

Using

\[
\begin{split}
{\rm Li}_3 \left[ { \pm \frac{1}{3}} \right] &= \frac{1}{{3^6 }}\sum\limits_{k = 0}^\infty  {\frac{1}{{3^{6k} }}\left[ { \pm \frac{{3^5 }}{{(6k + 1)^3 }} + \frac{{3^4 }}{{(6k + 2)^3 }} \pm \frac{{3^3 }}{{(6k + 3)^3 }}} \right.}\\
&\left. { + \frac{{3^2 }}{{(6k + 4)^3 }} \pm \frac{3}{{(6k + 5)^3 }} + \frac{1}{{(6k + 6)^3 }}} \right]
\end{split}
\]

leads to the ternary BBP-type formula

\[
\begin{split}
&13\zeta (3) - \pi ^2 \ln 3 + \ln ^3 3\\
&\quad= \frac{2}{{3^5 }}\sum\limits_{k = 0}^\infty  {\frac{1}{{3^{6k} }}\left[ {\frac{{3^6 }}{{(6k + 1)^3 }} + \frac{{3^4 }}{{(6k + 2)^3 }} + \frac{{3^4 }}{{(6k + 3)^3 }}} \right.}\\
&\qquad\left. { + \frac{{3^2 }}{{(6k + 4)^3 }} + \frac{{3^2 }}{{(6k + 5)^3 }} + \frac{1}{{(6k + 6)^3 }}} \right]\,.
\end{split}
\]

A shorter version (length~2) of the above formula is

\[
13\zeta (3) - \pi ^2 \ln 3 + \ln ^3 3 = \frac{2}{3}\sum\limits_{k = 0}^\infty  {\frac{1}{{9^k }}\left[ {\frac{9}{{(2k + 1)^3 }} + \frac{1}{{(2k + 2)^3 }}} \right]}\,. 
\]

The ternary BBP-type formula that results from~\eqref{equ.thgd6ju} is discussed elsewhere~\cite{adenyjm}. 

\bigskip
 
Next we obtain the ternary BBP-type formula that results from~\eqref{equ.hv0t3gq}.

\begin{equation}\label{equ.retygg}
\begin{split}
{\mathop{\rm Im\,}\nolimits} {\rm Li}_3 \left[ {\frac{1}{{\sqrt 3 }}\exp \left( {\frac{{i\pi }}{2}} \right)} \right] &= \frac{{\sqrt 3 }}{{3^6 }}\sum\limits_{k = 0}^\infty  {\frac{1}{{3^{6k} }}\left[ {\frac{{3^5 }}{{(12k + 1)^3 }} - \frac{{3^4 }}{{(12k + 3)^3 }}} \right.}\\
&\quad\left. { + \frac{{3^3}}{{(12k + 5)^3 }} - \frac{{3^2 }}{{(12k + 7)^3 }} + \frac{3}{{(12k + 9)^3 }} - \frac{1}{{(12k + 11)^3 }}} \right]
\end{split}
\end{equation}
and

\begin{equation}\label{equ.cfdrek}
\begin{split}
{\mathop{\rm Im\,}\nolimits} {\rm Li}_3 \left[ {\frac{1}{{\sqrt 3 }}\exp \left( {\frac{{i\pi }}{6}} \right)} \right] &= \frac{{\sqrt 3 }}{{2 \cdot 3^6 }}\sum\limits_{k = 0}^\infty  {\frac{1}{{3^{6k} }}\left[ {\frac{{3^5 }}{{(12k + 1)^3 }} + \frac{{3^5 }}{{(12k + 2)^3 }} + \frac{{2 \cdot 3^4 }}{{(12k + 3)^3 }}} \right.}\\
&+ \frac{{3^4 }}{{(12k + 4)^3 }} + \frac{{3^3 }}{{(12k + 5)^3 }} - \frac{{3^2 }}{{(12k + 7)^3 }}\\
&\left. { - \frac{{3^2 }}{{(12k + 8)^3 }} - \frac{{2 \cdot 3}}{{(12k + 9)^3 }} - \frac{3}{{(12k + 10)^3 }} - \frac{1}{{(12k + 11)^3 }}} \right]
\end{split}
\end{equation}

Combining~\eqref{equ.retygg} and \eqref{equ.cfdrek} according to the prescription of~\eqref{equ.hv0t3gq}, we arrive at

\[
\begin{split}
\frac{{29}}{{1296}}\pi ^3\sqrt 3- \frac{1}{{48}}\pi\sqrt 3 \ln ^2 3 &= \frac{{1 }}{{2 \cdot 3^5 }}\sum\limits_{k = 0}^\infty  {\frac{1}{{3^{6k} }}\left[ {\frac{{3^6 }}{{(12k + 1)^3 }} - \frac{{5 \cdot 3^5 }}{{(12k + 2)^3 }}} \right.}\\
&\qquad- \frac{{2 \cdot 3^6 }}{{(12k + 3)^3 }} - \frac{{5 \cdot 3^4 }}{{(12k + 4)^3 }} + \frac{{3^4 }}{{(12k + 5)^3 }} - \frac{{3^3 }}{{(12k + 7)^3 }}\\
&\qquad\qquad\left. { + \frac{{5 \cdot 3^2 }}{{(12k + 8)^3 }} + \frac{{2 \cdot 3^3 }}{{(12k + 9)^3 }} + \frac{{5 \cdot 3}}{{(12k + 10)^3 }} - \frac{3}{{(12k + 11)^3 }}} \right]\,.
\end{split}
\]

\section{Degree~$4$ Ternary BBP-type Formulas}

A two-variable functional equation for degree~4 polylogarithms (\mbox{equation A.2.7.40} of~\cite{lewin81}) reads 

\begin{equation}\label{equ.ieibdj8}
\begin{split}
&\operatorname{Li}_4 \left[-{\frac {{x}^{2}y\eta}{\xi}} \right] +\operatorname{Li}_4 \left[-{\frac {{
y}^{2}x\xi}{\eta}} \right] +\operatorname{Li}_4 \left[{\frac {{x}^{2}y}{{\eta}^{2}\xi
}} \right] +\operatorname{Li}_4 \left[{\frac {{y}^{2}x}{{\xi}^{2}\eta}} \right]\\ 
&=6\,\operatorname{Li}_4\left[xy \right] +6\,\operatorname{Li}_4 \left[{\frac {xy}{\eta\,\xi}} \right] +6
\,\operatorname{Li}_4 \left[-{\frac {xy}{\eta}} \right] +6\,\operatorname{Li}_4 \left[-{\frac {xy}{
\xi}} \right]\\
&+3\,\operatorname{Li}_4 \left[x\eta \right] +3\,\operatorname{Li}_4 \left[y\xi
 \right] +3\,\operatorname{Li}_4 \left[{\frac {x}{\eta}} \right] +3\,\operatorname{Li}_4 \left[{
\frac {y}{\xi}} \right] +3\,\operatorname{Li}_4 \left[-{\frac {x\eta}{\xi}} \right]\\
& +3\,\operatorname{Li}_4 \left[-{\frac {y\xi}{\eta}} \right] +3\,\operatorname{Li}_4 \left[-{\frac {x}
{\eta\,\xi}} \right] +3\,\operatorname{Li}_4 \left[-{\frac {y}{\eta\,\xi}} \right] -6
\,\operatorname{Li}_4 \left[x \right]\\
& -6\,\operatorname{Li}_4 \left[y \right] -6\,\operatorname{Li}_4 \left[-{
\frac {x}{\xi}} \right] -6\,\operatorname{Li}_4 \left[-{\frac {y}{\eta}} \right] +3/2
\,\ln^2\xi\ln^2\eta\,,
\end{split}
\end{equation}

where $\xi=1-x$, $\eta=1-y$.

\bigskip

Putting $x=-\exp {(i\pi/3)}$ and $y=\exp {(i\pi/3)}$ in~\eqref{equ.ieibdj8}, simplifying and taking real and imaginary parts, we obtain

\begin{equation}\label{equ.u7f27ai}
\begin{split}
&- 12\;{\mathop{\rm Re\;}\nolimits} {\rm Li}_4 \left[ {\frac{1}{{\sqrt 3 }}e^{i{\pi  \mathord{\left/
 {\vphantom {\pi  2}} \right.
 \kern-\nulldelimiterspace} 2}} } \right] - 3\;{\mathop{\rm Re\;}\nolimits} {\rm Li}_4 \left[ {\frac{1}{{\sqrt 3 }}e^{i{\pi  \mathord{\left/
 {\vphantom {\pi  6}} \right.
 \kern-\nulldelimiterspace} 6}} } \right]\\
 &\qquad+ {\rm Li}_4 \left[ {\frac{1}{3}} \right] + \frac{1}{4}\,{\rm Li}_4 \left[ { - \frac{1}{3}} \right]\\
&\qquad\qquad=  - \frac{{127\pi ^4 }}{{10368}} + \frac{1}{{64}}\pi ^2 \ln ^2 3 - \frac{5}{{384}}\ln ^4 3
\end{split}
\end{equation}

and

\begin{equation}\label{equ.xnlkxjb}
\begin{split}
&- 12\;{\mathop{\rm Im\;}\nolimits} {\rm Li}_4 \left[ {\frac{1}{{\sqrt 3 }}e^{i{\pi  \mathord{\left/
 {\vphantom {\pi  2}} \right.
 \kern-\nulldelimiterspace} 2}} } \right] + 15\;{\mathop{\rm Im\;}\nolimits} {\rm Li}_4 \left[ {\frac{1}{{\sqrt 3 }}e^{i{\pi  \mathord{\left/
 {\vphantom {\pi  6}} \right.
 \kern-\nulldelimiterspace} 6}} } \right]\\
&\qquad= \frac{{29}}{{864}}\pi ^3 \ln 3 - \frac{1}{{96}}\pi \ln ^3 3 - \frac{{11}}{3}{\rm Cl}_4 \left( {\frac{\pi }{3}} \right)\,.
\end{split}
\end{equation}

First we proceed to obtain the BBP-type formula invoked by~\eqref{equ.u7f27ai}.

\bigskip

Now,
\begin{equation}\label{equ.jrn4eaz}
\begin{split}
&{\mathop{\rm Re\,}\nolimits} {\rm Li}_4 \left[ {\frac{1}{{\sqrt 3 }}\exp \left( {\frac{{i\pi }}{2}} \right)} \right] = \frac{1}{{3^6 }}\sum\limits_{k = 0}^\infty  {\frac{1}{{3^{6k} }}\left[ { - \frac{{3^5 }}{{(12k + 2)^4 }} + \frac{{3^4 }}{{(12k + 4)^4 }}} \right.}\\
&\qquad\qquad\qquad\qquad\qquad\qquad\left. { - \frac{{3^3 }}{{(12k + 6)^4 }} + \frac{{3^2 }}{{(12k + 8)^4 }} - \frac{3}{{(12k + 10)^4 }} + \frac{1}{{(12k + 12)^4 }}} \right]
\end{split}
\end{equation}

and

\begin{equation}\label{equ.qpxxzi4}
\begin{split}
{\mathop{\rm Re\,}\nolimits} {\rm Li}_4 \left[ {\frac{1}{{\sqrt 3 }}\exp \left( {\frac{{i\pi }}{6}} \right)} \right]& = \frac{1}{{2 \cdot 3^6 }}\sum\limits_{k = 0}^\infty  {\frac{1}{{3^{6k} }}\left[ {\frac{{3^6 }}{{(12k + 1)^4 }} + \frac{{3^5 }}{{(12k + 2)^4 }}} \right.}\\
&\qquad- \frac{{3^4 }}{{(12k + 4)^4 }} - \frac{{3^4 }}{{(12k + 5)^4 }} - \frac{{2 \cdot 3^3 }}{{(12k + 6)^4 }}\\
&\qquad\qquad- \frac{{3^3 }}{{(12k + 7)^4 }} - \frac{{3^2 }}{{(12k + 8)^4 }} + \frac{3}{{(12k + 10)^4 }}\\
&\qquad\qquad\qquad\left. { + \frac{2}{{(12k + 12)^4 }}} \right]\,.
\end{split}
\end{equation}

Also

\begin{equation}\label{equ.jwslz5c}
\begin{split}
{\rm Li}_4 \left[ { \pm \frac{1}{3}} \right] &= \frac{{2^4 }}{{3^6 }}\sum\limits_{k = 0}^\infty  {\frac{1}{{3^{6k} }}\left[ { \pm \frac{{3^5 }}{{(12k + 2)^4 }} + \frac{{3^4 }}{{(12k + 4)^4 }}} \right.}\\
&\qquad\pm \frac{{3^3 }}{{(12k + 6)^4 }} + \frac{{3^2 }}{{(12k + 8)^4 }} \pm \frac{3}{{(12k + 10)^4 }}\\
&\qquad\qquad\left. { + \frac{1}{{(12k + 12)^4 }}} \right]\,.
\end{split}
\end{equation}

Combining~\eqref{equ.jrn4eaz}, \eqref{equ.qpxxzi4} and \eqref{equ.jwslz5c} according to~\eqref{equ.u7f27ai}, we obtain the following degree~4 ternary BBP-type formula:
\[
\begin{split}
\frac{{127\pi ^4 }}{{5184}} - \frac{{\pi ^2 \ln ^2 3}}{{32}} + \frac{{5\ln ^4 3}}{{192}} &= \frac{1}{{3^6 }}\sum\limits_{k = 0}^\infty  {\frac{1}{{3^{6k} }}\left[ {\frac{{3^7 }}{{(12k + 1)^4 }} - \frac{{5 \cdot 3^7 }}{{(12k + 2)^4 }}} \right.}\\
&\qquad - \frac{{19 \cdot 3^4 }}{{(12k + 4)^4 }} - \frac{{3^5 }}{{(12k + 5)^4 }} - \frac{{2 \cdot 3^6 }}{{(12k + 6)^4 }}\\
&\qquad\qquad - \frac{{3^4 }}{{(12k + 7)^4 }} - \frac{{19 \cdot 3^2 }}{{(12k + 8)^4 }} - \frac{{5 \cdot 3^3 }}{{(12k + 10)^4 }}\\
&\qquad\qquad\qquad\left. { + \frac{{3^2 }}{{(12k + 11)^4 }} - \frac{{10}}{{(12k + 12)^4 }}} \right]\,.
\end{split}
\]

Next we obtain the BBP-type formula invoked by~\eqref{equ.xnlkxjb}.

\bigskip

Writing

\begin{equation}\label{equ.zfwoayv}
\begin{split}
&{\mathop{\rm Im\,}\nolimits} {\rm Li}_4 \left[ {\frac{1}{{\sqrt 3 }}\exp \left( {\frac{{i\pi }}{2}} \right)} \right]\\
&\qquad = \frac{{\sqrt 3 }}{{27}}\sum\limits_{k = 0}^\infty  {\left( { - \frac{1}{{27}}} \right)^k \left[ {\frac{9}{{(6k + 1)^4 }} - \frac{3}{{(6k + 3)^4 }} + \frac{1}{{(6k + 5)^4 }}} \right]}
\end{split}
\end{equation}

and

\[
\begin{split}
{\mathop{\rm Im\,}\nolimits} {\rm Li}_4 \left[ {\frac{1}{{\sqrt 3 }}\exp \left( {\frac{{i\pi }}{6}} \right)} \right] &= \frac{{\sqrt 3 }}{{54}}\sum\limits_{k = 0}^\infty  {\left( { - \frac{1}{{27}}} \right)^k \left[ {\frac{9}{{(6k + 1)^4 }} + \frac{9}{{(6k + 2)^4 }}} \right.}\\
&\qquad\qquad\left. { + \frac{6}{{(6k + 3)^4 }} + \frac{3}{{(6k + 4)^4 }} + \frac{1}{{(6k + 5)^4 }}} \right]\,,
\end{split}
\]

and combining these according to~\eqref{equ.xnlkxjb}, we obtain the following degree~4 BBP-type formula:
\begin{equation}\label{equ.accm6ov}
\begin{split}
&\frac{1}{\sqrt 3}\left(11{\rm Cl}_4 \left( {\frac{\pi }{3}} \right) - \frac{{29}}{{288}}\pi ^3 \ln 3 + \frac{{\pi \ln ^3 3}}{{32}}\right)\\
&\qquad= \frac{{1}}{2}\sum\limits_{k = 0}^\infty  {\left( { - \frac{1}{{27}}} \right)^k \left[ {\frac{9}{{(6k + 1)^4 }} - \frac{{15}}{{(6k + 2)^4 }}} \right.}\\
&\qquad\qquad\left. { - \frac{{18}}{{(6k + 3)^4 }} - \frac{5}{{(6k + 4)^4 }} + \frac{1}{{(6k + 5)^4 }}} \right]\,.
\end{split}
\end{equation}

It is interesting to remark that~\eqref{equ.accm6ov} was also obtained by Broadhurst~\cite{broadhurst98a}, using the PSLQ Algorithm. We have thus found its formal proof for the first time, through~\eqref{equ.xnlkxjb}! 

\section{Degree~$5$ Ternary BBP-type Formulas}
The following degree~5 polylogarithm identity is derived in~\cite{broadhurst98}
\[
\begin{split}
\label{equ.t9v474l}
&{\rm Li}_5 \left[{\frac {x\alpha}{y\beta}} \right] +{\rm Li}_5 \left[x\alpha\,y
\eta \right] +{\rm Li}_5 \left[{\frac {x\alpha\,\beta}{\eta}} \right] +{\rm Li}_5
 \left[x\xi\,y\beta \right] +{\rm Li}_5 \left[{\frac {x\xi}{y\eta}}
 \right]\\
 &+{\rm Li}_5 \left[{\frac {x\xi\,\eta}{\beta}} \right] +{\rm Li}_5 \left[{
\frac {\alpha\,y\beta}{\xi}} \right] +{\rm Li}_5 \left[{\frac {\alpha}{\xi\,
y\eta}} \right] +{\rm Li}_5 \left[{\frac {\alpha\,\eta}{\xi\,\beta}}
 \right]\\
 &-9\,{\rm Li}_5 \left[xy \right] -9\,{\rm Li}_5 \left[x\beta \right] -9\,{\rm Li}_5
 \left[x\eta \right] -9\,{\rm Li}_5 \left[{\frac {x}{y}} \right] -9\,{\rm Li}_5
 \left[{\frac {x}{\beta}} \right]\\ 
&-9\,{\rm Li}_5 \left[{\frac {x}{\eta}}
 \right] -9\,{\rm Li}_5 \left[\alpha\,y \right] -9\,{\rm Li}_5 \left[\alpha\,\beta
 \right] -9\,{\rm Li}_5 \left[\alpha\,\eta \right]\\ 
&-9\,{\rm Li}_5 \left[{\frac {
\alpha}{y}} \right] -9\,{\rm Li}_5 \left[{\frac {\alpha}{\beta}} \right] -9
\,{\rm Li}_5 \left[{\frac {\alpha}{\eta}} \right] -9\,{\rm Li}_5 \left[\xi\,y
 \right] -9\,{\rm Li}_5 \left[\xi\,\beta \right]\\
 &-9\,{\rm Li}_5 \left[\xi\,\eta
 \right] -9\,{\rm Li}_5 \left[{\frac {y}{\xi}} \right] -9\,{\rm Li}_5 \left[{
\frac {\beta}{\xi}} \right] -9\,{\rm Li}_5 \left[{\frac {\eta}{\xi}}
 \right]\\
&  +18\,{\rm Li}_5 \left[x \right] +18\,{\rm Li}_5 \left[\alpha \right] +18
\,{\rm Li}_5 \left[\xi \right] +18\,{\rm Li}_5 \left[y \right] +18\,{\rm Li}_5 \left[
\beta \right]\\
&+18\,{\rm Li}_5 \left[\eta \right] -18\,\zeta (5) =
3/10\, \left( \ln\xi\right) ^{5}+3/4\, \left( \ln y -\ln x  \right)  \left( \ln\xi  \right) ^{4}\\
&+3/2\, \left( 3\,\ln y -\ln\eta  \right)  \left( \ln\eta\right) ^{2} \left( \ln \xi  \right) ^{2}+1/2\,{\pi }
^{2} \left( \ln \xi -3\,\ln \eta\right)\left( \ln \xi\right)^{2}+1/5\,{\pi }^{4}\ln  \xi\,.
\end{split}
\]

Here $\xi=1-x$, $\eta=1-y$, $\alpha=-x/\xi$ and $\beta=-y/\eta$.

\bigskip

Putting $x=-\exp {(i\pi/3)}$ and $y=\exp {(i\pi/3)}$ in~\eqref{equ.t9v474l} and simplifying, gives

\begin{equation}\label{equ.xmekn9g}
\begin{split}
&\frac{1}{{64}}\pi ^2 \ln ^3 3 - \frac{{127}}{{3456}}\pi ^4 \ln 3 - \frac{1}{{128}}\ln ^5 3 + \frac{{1573}}{{144}}\zeta (5)\\ 
  & = \frac{3}{2}\,{\rm Li}_5 \left[ { - \frac{1}{3}} \right] - 3\,{\rm Li}_5 \left[ {\frac{1}{3}} \right] + 9\,{\rm Li}_5 \left[ {\frac{1}{{\sqrt 3 }}e^{i{\pi  \mathord{\left/
 {\vphantom {\pi  6}} \right.
 \kern-\nulldelimiterspace} 6}} } \right] + 9\,{\rm Li}_5 \left[ {\frac{1}{{\sqrt 3 }}e^{ - i{\pi  \mathord{\left/
 {\vphantom {\pi  6}} \right.
 \kern-\nulldelimiterspace} 6}} } \right]\,.
\end{split}
\end{equation}

On taking real parts

\begin{equation}\label{equ.xmekn9g}
\begin{split}
&\frac{1}{{64}}\pi ^2 \ln ^3 3 - \frac{{127}}{{3456}}\pi ^4 \ln 3 - \frac{1}{{128}}\ln ^5 3 + \frac{{1573}}{{144}}\zeta (5)\\ 
&= 18\,{\mathop{\rm Re\;}\nolimits} {\rm Li}_5 \left[ {\frac{1}{{\sqrt 3 }}e^{i{\pi  \mathord{\left/
 {\vphantom {\pi  6}} \right.
 \kern-\nulldelimiterspace} 6}} } \right] + \frac{3}{2}\,{\rm Li}_5 \left[ { - \frac{1}{3}} \right] - 3\,{\rm Li}_5 \left[ {\frac{1}{3}} \right]\,.
\end{split}
\end{equation}

Now

\begin{equation}\label{equ.fxqzn9v}
\begin{split}
&18\,{\mathop{\rm Re\;}\nolimits} {\rm Li}_5 \left[ {\frac{1}{{\sqrt 3 }}e^{i{\pi  \mathord{\left/
 {\vphantom {\pi  6}} \right.
 \kern-\nulldelimiterspace} 6}} } \right] = 18\sum\limits_{k = 0}^\infty  {\left( {\frac{1}{{\sqrt 3 }}} \right)^k \frac{1}{{k^5 }}\cos } \left( {\frac{{k\pi }}{6}} \right)\\
&= \frac{1}{{3^6 }}\sum\limits_{k = 0}^\infty  {\frac{1}{{3^{6k} }}\left[ {\frac{{3^8 }}{{(12k + 1)^5 }} + \frac{{3^7 }}{{(12k + 2)^5 }} - \frac{{3^6 }}{{(12k + 4)^5 }}} \right.}\\
&- \frac{{3^6 }}{{(12k + 5)^5 }} - \frac{{2 \cdot 3^5 }}{{(12k + 6)^5 }} - \frac{{3^5 }}{{(12k + 7)^5 }} - \frac{{3^4 }}{{(12k + 8)^5 }}\\
&\left. { + \frac{{3^3 }}{{(12k + 10)^5 }} + \frac{{3^3 }}{{(12k + 11)^5 }} + \frac{{2 \cdot 3^2 }}{{(12k + 12)^5 }}} \right]\,,
\end{split}
\end{equation}

\begin{equation}\label{equ.u7hf4fg}
\begin{split}
&\frac{3}{2}\,{\rm Li}_5 \left[ { - \frac{1}{3}} \right] = \frac{3}{2}\sum\limits_{k = 1}^\infty  {\left( { - \frac{1}{3}} \right)^k \frac{1}{{k^5 }}}\\
&= \frac{1}{{3^6 }}\sum\limits_{k = 0}^\infty  {\frac{1}{{3^{6k} }}\left[ {\frac{{ - 3^6  \cdot 2^4 }}{{(12k + 2)^5 }} + \frac{{3^5  \cdot 2^4 }}{{(12k + 4)^5 }} - \frac{{3^4  \cdot 2^4 }}{{(12k + 6)^5 }}} \right.}\\
&\left. { + \frac{{3^3  \cdot 2^4 }}{{(12k + 8)^5 }} - \frac{{3^2  \cdot 2^4 }}{{(12k + 10)^5 }} + \frac{{3 \cdot 2^4 }}{{(12k + 12)^5 }}} \right]
\end{split}
\end{equation}

and

\begin{equation}\label{equ.s8q7u6v}
\begin{split}
&3\,{\rm Li}_5 \left[ {\frac{1}{3}} \right] = 3\sum\limits_{k = 1}^\infty  {\left( {\frac{1}{3}} \right)^k \frac{1}{{k^5 }}}\\
&= \frac{1}{{3^6 }}\sum\limits_{k = 0}^\infty  {\frac{1}{{3^{6k} }}\left[ {\frac{{3^6  \cdot 2^5 }}{{(12k + 2)^5 }} + \frac{{3^5  \cdot 2^5 }}{{(12k + 4)^5 }} + \frac{{3^4  \cdot 2^5 }}{{(12k + 6)^5 }}} \right.}\\
&\left. { + \frac{{3^3  \cdot 2^5 }}{{(12k + 8)^5 }} + \frac{{3^2  \cdot 2^5 }}{{(12k + 10)^5 }} + \frac{{3 \cdot 2^5 }}{{(12k + 12)^5 }}} \right]\,.
\end{split}
\end{equation}

Using~\eqref{equ.fxqzn9v}, \eqref{equ.u7hf4fg} and \eqref{equ.s8q7u6v} in~\eqref{equ.xmekn9g}, we obtain the ternary BBP-type formula
\[
\begin{split}
&\frac{1}{{64}}\pi ^2 \ln ^3 3 - \frac{{127}}{{3456}}\pi ^4 \ln 3 - \frac{1}{{128}}\ln ^5 3 + \frac{{1573}}{{144}}\zeta (5)\\
&= \frac{1}{{3^5 }}\sum\limits_{k = 0}^\infty  {\frac{1}{{3^{6k} }}\left[ {\frac{{3^7 }}{{(12k + 1)^5 }} - \frac{{5 \cdot 3^7 }}{{(12k + 2)^5 }} - \frac{{19 \cdot 3^4 }}{{(12k + 4)^5 }}} \right.}\\
&- \frac{{3^5 }}{{(12k + 5)^5 }} - \frac{{2 \cdot 3^6 }}{{(12k + 6)^5 }} - \frac{{3^4 }}{{(12k + 7)^5 }}\\
&\left. { - \frac{{19 \cdot 3^2 }}{{(12k + 8)^5 }} - \frac{{5 \cdot 3^3 }}{{(12k + 10)^5 }} + \frac{{3^2 }}{{(12k + 11)^5 }} - \frac{{10}}{{(12k + 12)^5 }}} \right]\,.
\end{split}
\]

\section{Conclusion}

Using a fairly straightforward method, we have obtained several ternary BBP-type formulas, which can now be added to the literature. In particular we proved the following formulas (written in the now standard BBP notation~\cite{bailey01}).\\
\\
$\ln 2=1/(2\cdot3^5){\rm P}((1,3^6,12,(3^5,3^5,0,-3^4,-3^3,0,-3^2,-3^2,0,3,1,0)))$\\
\\
$\ln 3=1/3^6{\rm P}((1,3^6,12,(3^6,3^5,0,-3^4,-3^4,-2\cdot 3^3,-3^3,-3^2,0,3,3,2)))$\\
\\
$\ln 3=1/27{\rm P}((1,-27,6,(27,9,0,-3,-3,-2)))$\\
\\
$\pi\sqrt 3=1/3^4{\rm P}((1,3^6,12,(3^5,3^5,2\cdot3^4,3^4,3^3,0,-3^2,-3^2,-6,-3,-1,0)))$\\
\\
$\pi\sqrt 3=2/3^4{\rm P}((1,3^6,12,(3^5,0,-3^4,0,3^3,0,-3^2,0,3,0,-1,0)))$\\
\\
$\pi\sqrt 3=6 {\rm P}((1,-3,2,(1,0)))$\\
\\
$\pi^2=2/27{\rm P}((2,3^6,12,(3^5,-5\cdot3^4,0,-3^4,-3^3,-2^3\cdot3^2,-3^2,-3^2,0,-5,1,0)))$\\
\\
$\ln^2 3=1/3^6{\rm P}((2,3^6,12,(2\cdot3^7,-2\cdot3^8,0,-2\cdot13\cdot3^4,-2\cdot3^5,-2^3\cdot3^5,-2\cdot3^4,\\
-2\cdot13\cdot3^2,0,-2\cdot3^4,2\cdot3^2,-8)))$\\
\\
$\pi\sqrt 3\ln 3=2{\rm P}((2,-27,6,(9,-15,-18,-5,1,0)))$\\
\\
$\sqrt 3{\rm Cl}_2(\pi/3)=1/3{\rm P}((2,-27,6,(9,-9,-12,-3,1,0)))$\\
\\
$13\zeta(3)-\pi^2\ln 3+\ln^3 3=2/3^5{\rm P}((3,3^6,6,(3^6,3^4,3^4,3^2,3^2,1)))$\\
\\
$13\zeta(3)-\pi^2\ln 3+\ln^3 3=2/3{\rm P}((3,9,2,(9,1)))$\\
\\
$(29\pi^3/1296-\pi\ln^2 3/48)/\sqrt 3=1/2/3^5{\rm P}((3,3^6,12,(3^6,-5\cdot3^5,\\
-2\cdot3^6,-5\cdot3^4,3^4,0,-3^3,5\cdot3^2,2\cdot3^3,5\cdot3,-3,0)))$\\
\\
$127\pi^4/5184-\pi^2\ln^2 3/32+5\ln^4 3/192=1/3^6{\rm P}((4,3^6,12,(3^7,-5\cdot3^7,0,\\
-19\cdot3^4,-3^5,-2\cdot3^6,-3^4,-19\cdot3^2,0,-5\cdot3^3,3^2,-10)))$\\
\\
$(11{\rm Cl}_4(\pi/3)-29\pi^3\ln 3/288+\pi\ln^3 3/32)/\sqrt 3=1/2{\rm P}((4,-27,6,(9,-15,-18,-5,1,0)))$\\
\\
$\pi^2\ln^3 3/64-127\pi^4\ln 3/3456-\ln^5 3/128+1573\zeta(5)/144=1/3^5{\rm P}((5,3^6,12,(3^7,\\
-5\cdot3^7,0,-19\cdot3^4,-3^5,-2\cdot3^6,-3^4,-19\cdot3^2,0,-5\cdot3^3,3^2,-10)))$

\bigskip

We also proved the following ternary zero relations:\\
\\
$0={\rm P}((1,3^6,12,(3^5,-5\cdot 3^4,0,-3^4,-3^3,-2^3\cdot 3^2,-3^2,-3^2,0,-5,1,0)))$\\
\\
$0={\rm P}((1,3^6,12,(3^5,-3^5,-2^2\cdot 3^4,-3^4,3^3,0,-3^2,3^2,3\cdot 2^2,3,-1,0)))$\\

\section*{Acknowledgments}
The author thanks Dr. D. H. Bailey, the foremost expert on BBP-type formulas and Experimental Mathematics, for interesting communications.

\end{document}